\documentclass[12pt, a4paper]{article}
\usepackage{amsfonts}
\usepackage{mathrsfs}
\usepackage{latexsym}
\usepackage{xy}
\usepackage{amsfonts,amsmath,amssymb,amsthm}
\usepackage{color}

\xyoption{all}

\newcommand{\bcen}{\begin{center}}     \newcommand{\ecen}{\end{center}}
\newcommand{\bay}{\begin{array}}      \newcommand{\eay}{\end{array}}
\newcommand{\beq}{\begin{eqnarray*}}      \newcommand{\eeq}{\end{eqnarray*}}

\def\cA{{\cal A}}

\def\cI{{\cal I}}
\def\cJ{{\cal J}}

\def\cT{{\cal T}}

\def\dim{\mathrm{dim}}
\def\mod{\mathrm{mod}}

\def\Im{\mathrm{Im}}
\def\Ker{\mathrm{Ker}}

\def\add{\mathrm{add}}

\begin{document}

\newtheorem{theorem}{Theorem}
\newtheorem{proposition}{Proposition}
\newtheorem{lemma}{Lemma}
\newtheorem{corollary}{Corollary}
\newtheorem{remark}{Remark}
\newtheorem{example}{Example}
\newtheorem{definition}{Definition}
\newtheorem{question}{Question}

\title{Derived dimensions of representation-finite algebras}

\author{Yang Han \footnote{The author is supported by NSFC (Project
10731070).}}

\date{\footnotesize KLMM, AMSS,
Chinese Academy of Sciences, Beijing 100190, China. \\ E-mail:
hany@iss.ac.cn}

\maketitle

\begin{abstract} It is shown that the derived dimension of any
representation-finite Artin algebra is at most one.
\end{abstract}


\bigskip

Let $\cT$ be a triangulated category, and $\cI$ and $\cJ$ two full
subcategories of $\cT$. Denote by $\langle \cI \rangle$ the smallest
full subcategory of $\cT$ containing $\cI$ and closed under shifts,
finite direct sums and direct summands. Denote by $\cI \ast \cJ$ the
full subcategory of $\cT$ consisting of all the objects $M \in \cT$
for which there exists a triangle $I \rightarrow M \rightarrow J
\rightarrow I[1]$ with $I \in \cI$ and $J \in \cJ$. Put $\cI
\diamond \cJ = \langle \cI \ast \cJ \rangle$ and inductively
$$\langle
\cI \rangle _n = \left\{\begin{array}{ll} \{0\},& \mbox{ if } n = 0;\\
\langle \cI \rangle _{n-1} \diamond \langle \cI \rangle , & \mbox{
if } n \geq 1.
\end{array}\right.$$
The {\it dimension} of $\cT$ is $\dim \cT := \mbox{min}\{d | \exists
M \in \cT, \ni \cT = \langle M \rangle _{d+1} \}$, or $\infty$ if
there is no such an $M$ for any $d$. (ref. \cite[Definition
3.2]{R}).

Let $\cA$ be an abelian category. Then $D^b(\cA)$, the derived
category of bounded complexes over $\cA$, is a triangulated
category. We call $\dim D^b(\cA)$ the {\it derived dimension} of
$\cA$. Let $A$ be an associative algebra with identity, and $A$-mod
the category of finitely generated left $A$-modules. Then we also
say $\dim D^b(A\!\!-\!\!\mod)$ is the {\it derived dimension} of
$A$. (ref. \cite{CYZ}). The derived dimension of an Artin algebra is
closed related to its Loewy length, global dimension, and
representation dimension. Especially, it provides a lower bound of
the representation dimension, which leads to find the examples of
the algebras of arbitrarily large representation dimensions. (ref.
\cite{R,KK}).

The following theorem is crucial though its proof is quite simple:

\medskip

\noindent{\bf Theorem.} {\it Let $\cA$ be an abelian category. Then
$D^b(\cA) = \langle \cA \rangle _2$. Furthermore, if $\cA = \add M$
for some object $M \in \cA$ then $\dim D^b(\cA) \leq 1$. Here $\add
M$ is the full subcategory of $\cA$ consisting of direct summands of
direct sums of finite copies of $M$.}

\begin{proof}
Let $X^\bullet = (X_n,f_n)_{n \in \mathbb{Z}}$ be a bounded complex
over $\cA$. Define the complexes $K^\bullet := (\Ker f_n,0)_{n \in
\mathbb{Z}}$ and $I^\bullet := (\Im f_n,0)_{n \in \mathbb{Z}}$. Then
$$0 \rightarrow K^\bullet \stackrel{(i_n)_{n \in
\mathbb{Z}}}{\rightarrow} X^\bullet \stackrel{(\tilde{f}_n)_{n \in
\mathbb{Z}}}{\rightarrow} I^\bullet \rightarrow 0$$ is an exact
triple of complexes, where $i_n : \Ker f_n \rightarrow X_n$ is the
natural embedding and $\tilde{f}_n : X_n \rightarrow \Im f_n$ is the
restriction of $f_n$ to $\Im f_n$ for all $n \in \mathbb{Z}$.
Namely, we have the following commutative diagram in which each
column is a short exact sequence:
$$\begin{array}{clclclc}
& 0 &  & 0 &  & 0 &  \\
& \downarrow & & \downarrow & & \downarrow & \\
\cdots \stackrel{0}{\rightarrow} & \Ker f_{n+1} &
\stackrel{0}{\rightarrow} & \Ker f_n & \stackrel{0}{\rightarrow} &
\Ker f_{n-1} & \stackrel{0}{\rightarrow} \cdots \\
& \downarrow i_{n+1} & & \downarrow i_n & &
\downarrow i_{n-1} & \\
\cdots \stackrel{f_{n+2}}{\rightarrow} & X_{n+1} &
\stackrel{f_{n+1}}{\rightarrow} & X_n & \stackrel{f_n}{\rightarrow}
& X_{n-1} & \stackrel{f_{n-1}}{\rightarrow} \cdots \\
& \downarrow \tilde{f}_{n+1} & & \downarrow \tilde{f}_n & &
\downarrow \tilde{f}_{n-1} & \\
\cdots \stackrel{0}{\rightarrow} & \Im f_{n+1} &
\stackrel{0}{\rightarrow} & \Im f_n & \stackrel{0}{\rightarrow} &
\Im f_{n-1} & \stackrel{0}{\rightarrow} \cdots \\
& \downarrow & & \downarrow & & \downarrow & \\
& 0 &  & 0 &  & 0 &
\end{array}$$
By \cite[Chapter III, \S 3, 5 Proposition]{GM}, we have a triangle
$K^\bullet \rightarrow X^\bullet \rightarrow I^\bullet \rightarrow
K^\bullet[1]$ in $D^b(\cA)$. Since $X^\bullet$ is a bounded complex
over $\cA$, $K^\bullet = \oplus_{n \in \mathbb{Z}} (\Ker f_n)[n]$
and $I^\bullet = \oplus_{n \in \mathbb{Z}} (\Im f_n)[n]$ have only
finitely many nonzero summands. Hence $K^\bullet$ and $I^\bullet \in
\langle \cA \rangle$, and thus $X^\bullet \in \langle \cA \rangle
_2$.

Furthermore, if $\cA = \add M$ for some object $M \in \cA$ then
$D^b(\cA) = \langle M \rangle _2$. Thus $\dim D^b(\cA) \leq 1$.
\end{proof}

An Artin algebra $A$ is said to be {\it representation-finite} or
{\it of finite representation type} if up to isomorphism there is
only a finite number of indecomposable $A$-modules in $A$-mod. (ref.
\cite[p.120]{A}).

\medskip

\noindent{\bf Corollary.} {\it Let $A$ be a representation-finite
Artin algebra. Then $\dim D^b(A-\!\!\mod) \leq 1$.}

\begin{proof} Let $M$ be a direct sum of a complete set of
representatives of finitely generated indecomposable left
$A$-modules. Then $A\!\!-\!\!\mod = \add M$. By Theorem, we have
$\dim D^b(A\!\!-\!\!\mod) \leq 1$.
\end{proof}

\begin{remark} In \cite{O}, Oppermann posed an open question ---
Are there nonsemisimple algebras $A$ such that the equality holds in
the inequality $\mbox{\rm rep.dim} A$ $\geq \dim
D^b(A\!\!-\!\!\mod)$? For representation-finite nonsemisimple Artin
algebras, by Corollary and \cite[p.139,Proposition]{A}, we always
have $\dim D^b(A\!\!-\!\!\mod) \leq 1 < 2 = \mbox{\rm rep.dim} A$.
\end{remark}

A finite-dimensional algebra $A$ is said to be {\it derived finite}
if up to shift and isomorphism there is only a finite number of
indecomposable objects in $D^b(A\!-\!\mod)$. (cf. \cite[Definition
1]{BM}). Clearly, $\dim D^b(A\!-\!\mod)=0$ if and only if $A$ is
derived finite.

\begin{remark} In \cite[Remark 7.29]{R}, Rouquier posed a question ---
Which finite-dimensional algebras can have a derived category of
dimension 1? So far we have known that the hereditary algebras (ref.
\cite[Proposition 7.27]{R} or \cite[Proposition 2.6]{KK}), the
radical square zero algebras (ref. \cite[Proposition 7.37]{R}), and
the representation-finite algebras, which are not derived finite,
are of derived dimension one. \end{remark}

From now on we always consider finite-dimensional algebras over an
algebraically closed field!

\begin{remark} By Corollary and \cite[Theorem]{CYZ},
we can know explicitly the derived dimension of any
representation-finite algebra.
\end{remark}

\begin{remark} According to the derived category of bounded
complexes of left modules, there are four important classes of
algebras: derived finite algebras, derived discrete algebras,
derived tame algebras and derived wild algebras (ref.
\cite{BM,D,GK,V}). We have known that the derived finite algebras
are just the algebras of derived dimension zero. Let $A$ be a
derived discrete algebra. Then $\dim D^b(A\!-\!\mod) \leq 1$.
Indeed, by \cite[Theorem]{V}, $A$ is either a derived hereditary
algebra of Dynkin type or a representation-finite gentle algebra. In
both cases, $A$ is representation-finite. By Theorem, we have $\dim
D^b(A\!-\!\mod) \leq 1$. It is unknown if there is a upper bound for
the derived dimensions of derived tame or tame algebras.
\end{remark}

\end{document}